\documentclass[12pt]{article}
\usepackage[T2A]{fontenc}
\usepackage{amsmath}
\usepackage[cp1251]{inputenc}
\begin{document}
\large

\begin{center}
{\bf IDENTIFICATION OF BOUNDARY CONDITIONS USING NATURAL FREQUENCIES
IN CASE OF A RING MEMBRANE}
\end{center}

\vspace{0.2cm}

\centerline {Akhtymov A. M.$^{a}$, Mouftakhov A. V.$^{b}$\, Teicher
M.$^{b}$ }

\vspace{0.2cm}

$^a$ Institute of mechanics, Ufa, Russia;

$^b$ Department of Mathematics, Bar-Ilan University, Ramat-Gan,
Israel.

\vspace{0.2cm}

\date{}

\addtocounter{section}{1}
\begin{abstract}
The problem of finding boundary conditions for fastening of a ring
membrane, which are inaccessible for direct observation from the
natural frequencies of its flexural oscillations, is considered.

Two theorems on the uniqueness of this problem are proved, and a
method for establishing the unknown conditions for fastening of the
membrane to the walls is indicated. An approximate formula for
determining the unknown conditions is obtained, using
 first three natural frequencies.
 The method of approximate calculation of unknown boundary conditions, is
explained with the help of an example.

\vspace{0.2cm}


\end{abstract}
\vspace{0.2cm}

\vspace{0.2cm} {\bf 1. Formulation of the Inverse Problem.}

 Let us
give the mathematical formulation of the problem, before describing
the method of its solution. The problem of the radial oscillations
of membrane, is reduced after making a separation of variables,  to
the following eigenvalue problem (see [5])
$$
\ddot{y}(r)+\dot{y}(r)/r +\lambda ^2\, y(r) = 0,\qquad\eqno(1)
$$
$$
U_1(y)=k_1\,\dot{y}(a)-k_2\, y(a) = 0,\,\,\,
U_2(y)=k_3\,\dot{y}(b)+k_4\, y(b) = 0,\qquad\eqno(2)
$$
where $\lambda$, is the spectral parameter, $a$ and $b$ are the
short and long radius of the ring membrane respectively.

Now we formulate the inverse of the eigenvalue problem:

{\bf Problem 1}(inverse problem).
 {\it It is required to find the unknown linear forms $U_1(y)$ and  $U_2(y)$ from
 the eigenvalues of the problem (1)--(2)}.

Let $L(y,\lambda)=\ddot{y}(r)+\dot{y}(r)/r +\lambda ^2\, y(r), \,\,
(U_1(y),U_2(y))^T=A\,X^T,\\$ where $A=\{a_{ij}| i=1,2;j=1,2,3,4\}$
is a matrix of order $2\times4$ such that $a_{11}=k_1, \,
a_{12}=-k_2, \,
 a_{23}=k_3, \, a_{24}=k_4, \, a_{13}=a_{14}=a_{21}=a_{22}=0$, and

$ X(y)=(\dot{y}(a),y(a),\dot{y}(b),y(b)) $ is a row-vector. Then the
eigenvalue problem (1)--(2) is equivalent to the following
eigenvalue problem
$$
L(y,\lambda)=0,\,\,\,A\,X^T(y)=0.\qquad\eqno(3)
$$
We consider also another eigenvalue problem
$$
L(y,\lambda)=0,\,\,\,\widetilde{A}\,X^T(y)=0.\qquad\eqno(4)
$$
In these problems only the matrices of boundary conditions are
different.

 {\bf Proposition 1.} {\it If matrix $\widetilde{A}$ is
 equivalent to matrix ${A}$, i.e. there exist non-singular matrix
$S$ such that $A=S \widetilde {A}$, then problem (4) is equivalent
to problem (3).}

The proof is trivial.

 {\bf Proposition 2.} {\it Matrix $\widetilde{A}$ is
 equivalent to matrix ${A}$, if and only if, their  correspond minors of
order 2 are equal, accurate to coefficient.}

The proof is in [3].

In the projective geometry we define the Pl\"ucker coordinates of
the set of equivalent matrices of order $m\times n$ to be any one of
the equivalent $C^{(m+1)}_{(n+1)} $-tuples of its minors of order
$m$.

Now we shell give the following

 {\bf Definition 1.} {\it Pl\"ucker coordinates } of boundary conditions of the eigenvalue problem (3) are called
 any one of the equivalent $8$-tuples of  minors (of order $2$) of its matrix
 $A$: $(\ldots ,\ A_{ij} \,,\, \ldots)$.

From Proposition 1 and Proposition 2 we get

{\bf Proposition 3.} {\it If Pl\"ucker coordinates  of boundary
conditions of the eigenvalue problems (3) and (4)  are identical,
then problems (3) and (4) are also identical.}

Then  Problem 1 is equivalent to the following problem:

 {\bf Problem 2}.
 {\it It is required to find Pl\"ucker coordinates of boundary
 conditions of the problem (3) from  its eigenvalues, if its differential equation is known}.

\vspace{0.2cm}

{\bf 2. The Uniqueness of the Solution of the Inverse Problem.}

 Together with problem ~(3), let us consider problem ~(4).

\vspace{0.2 cm}

{\bf Theorem 1 }(on the uniqueness of the solution of the inverse
problem). {\it Suppose the following conditions are satisfied:
$$
  {\rm rank\, } A={\rm rank\, } \tilde{A} =2.\qquad\eqno(5)
$$
  If non-zero eigenvalues of problems
(3) and (4) are identical,  with account taken for their
multiplicities,  then problems (3) and (4) are also identical. }

\vspace{0.1 cm}

{\bf Proof.} A general solution of the differential equation of the
problem (3) is the function $ y(r)= y(r,\lambda ) = C_1\,
J_0(\lambda\, r)+C_2\, Y_0(\lambda\, r), $ where standard notations
for cylinder functions are used (see [2]).

The boundary conditions  are used to determine the constants $C_1$
and $C_2$.

The following function is a characteristic determinant of the
problem (3):
$$ \Delta (\lambda ) \equiv \left|
\begin{array}{cc}
U_1 \Big( J_0(\lambda\, a)\Big) & U_1 \Big( Y_0(\lambda\, a)\Big)
\\
U_2\Big( J_0(\lambda\, b)\Big) & U_2\Big( Y_0(\lambda\, b)\Big)
\end{array}
\right| \label{IP PMM U1 Delta pl}
$$

$\Delta (\lambda ) \equiv \det(A\,B^T)$, where
 $B$ is a matrix of order $2\times4$ such that
$b_{11}=\dot{J}_0(\lambda\,a) \,$, $\,b_{12}=J_0(\lambda\,a)
\,$,$\,b_{13}=\dot{J}_0(\lambda\,b) \,$,$\,b_{14}=J_0(\lambda\,a)
\,$,$\,b_{21}= \dot{Y}_0(\lambda\,a)\,$,$\,b_{22}=Y_0(\lambda\,a)
\,$,$\,b_{23}=\dot{Y}_0(\lambda\, b) \,$,$\,b_{24}=Y_0(\lambda\, b)
\,$. The  equation for frequencies, is obtained from the condition
of the existence of a non-zero solution for $C_i.$. The latter
solution exists if and only if, the determinant  $\Delta (\lambda )$
is equal to zero (see [5]).

Using the Binet-Cauchy formula, we obtain
$$
\Delta (\lambda ) = \sum_{1\leq i<j\leq 4}\; A_{ij} \cdot
B_{ij}(\lambda) \equiv 0. \eqno{(6)}
$$
It follows from asymptotic estimations for cylindrical functions
$$
J_{\nu}(z)\sim \sqrt{\frac{2}{\pi\, z}}\, \cos (z-\nu\,\pi /2-\pi
/4), \quad Y_{\nu}(z)\sim \sqrt{\frac{2}{\pi\, z}}\, \sin
(z-\nu\,\pi /2-\pi /4),
$$
that $\Delta(\lambda ) $ is an entire function (see [4]). Hence, it
follows from Hadamard's factorization theorem (see [4]), that
characteristic determinant $\Delta(\lambda ) $ of problem (3) and
characteristic determinant $\widetilde  \Delta(\lambda ) $ of
problem (4), are connected by the relation
$$
\Delta (\lambda ) \equiv C\,\lambda ^{k}\, \widetilde \Delta(\lambda
), \eqno{(7)}
$$
where $k$ is a certain non-negative integer and $C$ is a certain
non-zero constant.
 From this, we obtain the identity
$$
 \sum_{1\leq i<j\leq 4}\; A_{ij} \cdot
B_{ij}(\lambda) \equiv C\,\lambda ^{k}\, \sum_{1\leq i<j\leq 4}\;
\widetilde A_{ij} \cdot B_{ij}(\lambda).
$$
It follows from this, that
$$
 \sum_{1\leq i<j\leq 4}\; (A_{ij} \cdot B_{ij}(\lambda) - \widetilde
 A_{ij} \cdot C\,\lambda ^{k}\,B_{ij}(\lambda))\equiv 0.
 \eqno{(8)}
$$

 Note, that  number $k$ in this identity is equal to zero.
Actually, let us assume the opposite: $k\ne 0$. Using Maple, we get
that the functions $B_{ij}(\lambda)$ $(1\leq i<j\leq 4)$ and also
the same functions multiplied by $\lambda ^{k}$, are linearly
independent.

From this and identity (8), we obtain $A_{ij}=\widetilde
 A_{ij}=0 \, (1\leq i<j\leq 4)$,
which contradict the condition (5) of the theorem.

Hence, $k=0$ and
$$
 \sum_{1\leq i<j\leq 4}\; (A_{ij}-C\,\widetilde
 A_{ij}) \cdot
B_{ij}(\lambda)\equiv 0.
 \eqno{(9)}
$$

 From this, by virtue of the
linearly independence of the corresponding functions, we obtain
$A_{ij}= C\,\widetilde  A_{ij}$, i.e. Pl\"ucker coordinates  of
boundary conditions of the eigenvalue problems (3) and (4)  are
identical.

Then, it follows from Proposition 3, that problems (3) and (4) are
also identical. The theorem is proved.

{\bf 3. Exact Solution of The Inverse Problem and Its Stability.}

This section is concerned with solving this problem and constructing
an exact solution of the inverse problem.

Suppose $\lambda_1$, $\lambda_2$, $\lambda_3$, are the eigenvalues
of problem (3). We substitute the values $\lambda_i$, $(i=1,2,3)$
into (6) . Since $A_{12}=A_{34}=0$,  we obtain a system of three
homogeneous algebraic equations in the four unknowns $A_{13}$,
$A_{14}$, $A_{23}$, $A_{24}$:
$$
F \cdot Z^T=0,  \eqno{(10)}
$$
where $F=\{f_{ij}| i=1,2;j=1,2,3,4\}$ is a matrix of order
$3\times4$ such that $f_{i1}=B_{13}(\lambda_i), \,
f_{12}=B_{14}(\lambda_i), \, f_{13}=B_{23}(\lambda_i), \,
a_{24}=B_{24}(\lambda_i), \, (i=1,2,3)$, and $
Z=(A_{13},A_{14},A_{23}, A_{24}$))  is a row-vector.

The resulting system has an infinite set of solutions. It follows
from the uniqueness theorem, which has been proved that the unknown
minors $A_{13}$, $A_{14}$, $A_{23}$, $A_{24}$ can be found apart
from a constant. Hence, the resulting system must have a rank of 3
and a solution, determined apart from a constant multiplier.

If the minors (Pl\"ucker coordinates)  $A_{13}$, $A_{14}$, $A_{23}$,
$A_{24}$ are found apart from a constant, then then Problem 2 is
solved. We can reconstruct any one of the equivalent matrix of
boundary conditions by its Pl\"ucker coordinates (see [3]).

For example, if $A_{13}=1$, $A_{12}=A_{34}=0$ then
 $$ A= \left\|\begin{array}{cccc} 1 & A_{23} & 0 & 0
\\
0 & 0 & 1 & A_{14}
\end{array}\right\| .
\eqno{(11)}
$$

This reasoning proves

{\bf Theorem 2.} {\it If the matrix $F$ of system (10) has a rank of
3, then the solution of the inverse problem of the reconstruction
boundary conditions~(3), is unique.}

Note, that theorem 2 is stronger than theorem 1. Theorem 2 used only
three natural frequencies for the reconstruction of the boundary
conditions and  not all natural frequencies as theorem 1, are used.

It is significant that we have the following theorem on stability of
the solution:

{\bf Theorem 3.} {\it Suppose that one of the third-order minors of
matrix $F$ is substantially non-zero. If $
|\widetilde{\lambda}_i-{\lambda}_i|<\delta<<1,$ then the boundary
conditions
 of problem (3), are close to
 the boundary conditions
 of problem (4).}

The proof is analogous to the proof of  the theorem on stability in
[1].

{\bf 4. Approximate Solution.}

Since small errors are possible when measuring natural frequencies,
the problem arises to find an algorithm for the approximate
determination of the type of fastening from the natural frequencies,
which are found with a certain error.

Let $ A^{\circ}_{13}$, $  A^{\circ}_{14} $, $ A^{\circ}_{23}$, $ A^{\circ}_{24}$ be
solution of the system of homogeneous algebraic equations(10). It is
unnecessary for these liquids to be totally immiscible.

Aside from measurement and calculation errors, it is unnecessary for
the values $  A^{\circ}_{13} $, $A^{\circ}_{14} $, $A^{\circ}_{23}$, $ A^{\circ}_{24}$
 to be minors of a matrix. So this problem is not trivial.
We must find minors  $  A_{13}$, $ A_{14} $, $A_{23}$, $ A_{24}$,
closer to the values  $  A^{\circ}_{13} $, $A^{\circ}_{14} $,
$A^{\circ}_{23}$, $ A^{\circ}_{24}$.

It is known from algebraic geometry [3] that the numbers $A_{12}$,
$A_{13}$, $A_{14}$, $A_{23}$, $A_{24}$, $A_{34}$ are minors of some
$2 \times 4$ matrix, if and only if, the following condition is
satisfied: $\, A_{12} A_{34} - A_{13} A_{24} + A_{14}A_{23} = 0.$
This condition is called the Pl\"ucker relation.

If the Plucker relation for these numbers is realized, then
 $   A^{\circ}_{13} $, $A^{\circ}_{14} $, $A^{\circ}_{23}$, $ A^{\circ}_{24}$
are minors of some matrix and corresponding boundary conditions are
found.

If the Pl\"ucker relation for numbers $ A^{\circ}_{12}$, $
A^{\circ}_{13} $, $A^{\circ}_{14} $, $A^{\circ}_{23}$, $
A^{\circ}_{24}$, $ A^{\circ}_{34}$
 is not realized, the required minors
$A_{12}$, $A_{13}$, $A_{14}$, $A_{23}$, $A_{24}$, $A_{34}$ are found
with the help of orthogonal projection. By definition, put $ \, x_1
= A_{12},
 \, x_2 = A_{34},
 \, x_3 = A_{13},
 \, x_4 = -A_{24},
\, x_5 = A_{14},
 \, x_6 = A_{23}.
$ Using this definition, we get the Pl\"ucker relation
$$
 x_1\, x_2 + x_3\, x_4+ x_5\, x_6  = 0, \eqno{(12)}
$$
which characterizes a surface $S$ in the 6-dimensional space. By
definition, put $ \, y_1 = A^{\circ}_{12},
 \, y_2 = A^{\circ}_{34},
 \, y_3 = A^{\circ}_{13},
 \, y_4 = -A^{\circ}_{24},
\, y_5 = A^{\circ}_{14},
 \, y_6 = A^{\circ}_{23}.
$

By definition, put
$$
\begin{array}{c}
X = (x_1, x_2 ,  x_3, x_4,  x_5, x_6) ,\quad Y = (y_1, y_2, y_3,
y_4, y_5, y_6),
\\
X^*=(x_2, x_1, x_4, x_3,  x_6, x_5),\quad Y^*=(y_2, y_1, y_4, y_3,
y_6, y_5),
\\
(\vec X,\vec Y)=x_1\, y_1 + x_2\, y_2+x_3\, y_3 + x_4\, y_4+x_5\,
y_5 + x_6\, y_6, \end{array}
$$
where $\quad \vec X=\overrightarrow{OX}$,  $\vec
Y=\overrightarrow{OY}$, $O$ is the origin.

Let $X$ be an orthogonal projection of $Y$ on surface (12). The
vector $\vec X^*$ is normal for surface (12) at the point $X$.
 It is identical to the following equations
$$
\vec Y = \vec X + p\,\vec X^*,\eqno{(13)} $$
 where $p$ is a real number.
Having solved a set of linear equations (13), with the unknowns
$x_1, x_2 ,  x_3, x_4,  x_5, x_6$, we obtain
$$
\vec X = \frac1{1-p^2}\, (\vec Y - p\, \vec Y^*).\eqno{(14)}
$$
From (15), it is easy to obtain (see [1])

$$p= \frac {(\vec Y,\vec Y)-\sqrt{(\vec Y,\vec Y)^2 -
(\vec Y,\vec Y^*)^2}}{(\vec Y,\vec Y ^*)}.\eqno{(15)}
$$

\vspace{0.2cm}

 {\bf 5. Example.}

If $\lambda_1= 2.93$, $\lambda_2= 6.16$, $\lambda_3= 9.34$
 corresponding to
the first three natural frequencies $\omega_i$
 determined, using instruments for measuring the natural frequencies,
  then the solution of system
 (10), apart from a constant, has the form
$
 A^{\circ}_{13} = 1.00\, C, \, A^{\circ}_{14} = 2.00\, C, \,
    A^{\circ}_{24}=-2.00\, C,\, A^{\circ}_{23} = -C. ,\
$ Using (14) and (15), we get $ A_{13}  =  0.72\, C,\, A_{14}  =
1.17\, C,\, A_{23}  =  -1.17\, C, \, A_{24}  =  -1.89\, C. \,$
Suppose $C=1/{A_{13}}$; then from (11), we obtain
$$
 A= \left\|\begin{array}{cccc} 1 & -1.62 & 0 & 0
\\
0 & 0 & 1 & 1.62
\end{array}\right\| .
$$

Note that the numbers $\lambda_1= 2.93$, $\lambda_2= 6.16$,
$\lambda_3= 9.34$ presented above are almost the same as the first
three exact values, that correspond to a case with
$$
 A= \left\|\begin{array}{cccc} 1 & -2 & 0 & 0
\\
0 & 0 & 1 & 2
\end{array}\right\| .
$$

 This means that the unknown
fixing of the membrane, which is inaccessible for direct
observation, has been satisfactorily determined. The errors of
determination of boundary conditions, are stipulated with round-off
errors of calculations.

 {\bf 6. Acknowledgements.}

This research was partially supported by the Russian Foundation for
Basic Research (06-01-00354a),  Emmy Noether Research Institute for
Mathematics, the Minerva Foundation of Germany, the Excellency
Center "Group Theoretic Methods in the Study of Algebraic Varieties"
 of the Israel Science Foundation, and by EAGER
(European Network in Algebraic Geometry).


\begin{thebibliography}{ABR}


\vspace{-0.3cm}
\bibitem{Akhtyamov 03 IPE}
{ A. \ M. \ Akhtyamov and A. \ V. \ Mouftakhov}, {\em Identification
of Boundary Conditions Using Natural Frequencies,} Inverse Probl. \
Sci. \ Eng.~{\bf 12(4)} (2004),393–408. 12(4), 393--408 (2004).

\vspace{-0.3cm}
\bibitem{Gontkevich 64}
{V. \ S. \ Gontkevich}, {\it Natural oscillations of plates and
shells}, “Naukova Dumka”, Kiev, 1964 (in Russian); German transl.,
VEB Fachbuchverlag, Leipzig, 1967.


\vspace{-0.3cm}
\bibitem{Hodge 94}
{ W. \ V. \ D. \ Hodge and D. \ Pedoe}, {\it Methods of Algebraic
Geometry}, Cambridge Univ. Press., Cambridge, UK, 1994.


\vspace{-0.3cm}

\bibitem{Gol'dberg 91}
{ B. \ Ya. \ Levin}, {\it Distribution of zeros of entire
functions}, Gostekhizdat,  Moscow, 1956. p. 632; English transl.:
Amer. Math. Soc., Providence, R.~I., 1980, p.~524.

\vspace{-0.3cm}

\bibitem{Rayleigh 29}
{ W. \ V. \ Strutt (Lord Rayleigh)}, {\it The Theory of Sound,} {\bf
vol.I}, Macmillan, London, 1929.

\end{thebibliography}
\end{document}